\documentclass[11pt,letterpaper]{article}
\usepackage[margin=1in]{geometry}
\date{}
\usepackage{color}
\usepackage{cite}
\usepackage{graphicx}
\usepackage[cmex10]{amsmath}
\usepackage{amssymb}
\usepackage{amsfonts}
\usepackage{float}
\usepackage{subfigure}
\usepackage{multirow}
\usepackage[capitalize]{cleveref}
\usepackage{algorithm}
\usepackage{algpseudocode}% http://ctan.org/pkg/algorithmicx
\usepackage{multirow}
\usepackage{multicol}
\usepackage{todonotes}
\usepackage{bm}

\usepackage{amsmath}

\usepackage{tikz}
\tikzstyle{cnode} = [draw, circle,scale=.6]
\tikzstyle{level 1} = [level distance=.35\textwidth, sibling distance=1\textwidth]
\tikzstyle{level 2} = [level distance=.35\textwidth, sibling distance=.45\textwidth]
\tikzstyle{level 3} = [level distance=.3\textwidth, sibling distance=.25\textwidth]
\usetikzlibrary{decorations.pathreplacing,patterns}

\usepackage[caption=false,font=footnotesize]{subfig}

\newcommand{\rev}[1]{{\color{black}{#1}}}

\begin{document}
	\title{Comb Tensor Network Finite-Difference Time-Domain Algorithms for 2D and 3D Electromagnetic Simulations}
	\author{Jonathan Tabares, Erika Ye, Yang Liu
		\thanks{\textit{Corresponding author: Yang Liu.}}
		\thanks{Y. Liu, and E. Ye are with Applied Mathematics and Computational Research Division, Lawrence Berkeley National Laboratory, Berkeley, CA 94720, USA. (e-mail: liuyangzhuan@lbl.gov, erikaye@lbl.gov). }
		\thanks{J. Tabares is with the Department of Physics, Florida International University, Miami, FL, 33199, USA. (e-mail: jtaba016@fiu.edu).}
	}
	
	\maketitle
	
	\begin{abstract}
		This paper presents efficient 2D and 3D finite-difference time-domain (FDTD) algorithms augmented by the comb tensor network (CTN) representation of the field and auxiliary quantities. The proposed CTN-FDTD algorithm leverages \rev{the compressibility of CTN} and rank-reduction strategies including local and global smoothing techniques to achieve\rev{, as long as the ranks remain bounded,} logarithmic memory requirement and computation (per time step) cost. Numerical examples validate these complexities \rev{when the ranks stay bounded, show where the ranks grow,} and demonstrate the applicability of CTN-FDTD using a variety of source excitation, frequency and permittivity profiles in both 2D and 3D. \rev{For a 3D dielectric sphere, the scattered field of CTN-FDTD with a maximum rank of 96 differs from that of standard FDTD by 2--4\%, and a quantized tensor train (QTT) chain run with the same implementation needs about three times the time per step to reach comparable accuracy.} Remarkably, CTN-FDTD can simulate 35 trillion spatial grid points in 3D with low memory and computation costs on one computation node, whereas standard FDTD would require 281 TB memory to store each of its 31 3D arrays. Moreover, CTN-FDTD was applied to realistic 2D photonic crystal structure simulation with one million grid points and 3D scattering analysis with \rev{68.7} billion grid points.

	\end{abstract}
	
	\par

	% \fi

	\section{Introduction}\label{sec:intro}

Finite-difference time-domain (FDTD) methods are among the most widely used numerical schemes for transient electromagnetic (EM) analysis because of their explicit time-marching structure, favorable scalability on structured grids, and ability to accommodate complex material profiles and absorbing boundary conditions \cite{atef2009finite,chew2001fast}. This simplicity has made FDTD a standard tool for wave propagation, scattering, and radiation analysis in optics, microwave engineering, and photonic design, etc. In the standard discretization, the computational cost and memory requirement scale as $O(N_tN^d)$ and $O(N^d)$, respectively, where $N$ is the number of grid points per dimension, $N_t$ is the number of time steps, and $d$ is the spatial dimensionality. However, large-scale FDTD simulations are limited by the well-known tradeoff between accuracy and cost: fine discretizations and dense representations of the perfectly matched layers (PML) or spatially-varying material parameters can suppress numerical dispersion and resolve geometry details, but lead to rapidly growing memory and CPU demands.

This challenge has motivated substantial interest in tensor-compressed EM solvers. In recent years, low-rank tensor networks have been used to accelerate a variety of differential and integral-equation formulations for large-scale electromagnetic simulations \cite{chen2019tensor,ye2026practical,ye2024quantized,jia2026id,zhou2025tensor,manzini2023tensor,giannakopoulos2022compression,kielstra2025linear,yucel2019new,chen2019sparsity}. Within FDTD-type solvers, tensor-train (TT) and quantized TT (QTT) representations have shown promising reductions in complexity for structured meshes and smooth solutions \cite{nguyen2026tensorFDTD,zhou2025tensor,manzini2023tensor,ye2026practical}. These methods exploit the fact that field quantities and material coefficients often exhibit low-rank structure in their multidimensional representation. Nonetheless, the effectiveness of such compression is highly sensitive to the ordering of tensor sites after quantization\rev{, and it depends on the compressibility of the field patterns. Zhou and Teixeira formulated 3D FDTD with PML in the TT format, performed the updates with rounding, and analyzed the rank growth caused by numerical noise \cite{zhou2025tensor}; the storage and the cost grow linearly with $N$. Nguyen and Okhmatovski incorporated QTT into 2D \cite{nguyen2026tensorFDTD} and full-vector 3D \cite{nguyen2026ttfdtd3d} scattered-field FDTD, with a logarithmic cost of the spatial operators when the ranks are bounded; they smooth the material interfaces with a Gaussian filter evaluated by a 3D FFT of the material array on the full grid, before the compression. }

This observation motivates the use of comb tensor networks (CTNs), which provide an alternative tensorization and quantization that are more naturally aligned with the geometry of multidimensional grid-based operators \cite{chepiga2019comb}. CTN organizes tensor indices with a structured ordering that preserves local spatial connectivity while retaining the compression advantages of tensor networks. In 2D, CTN is equivalent to QTT with a specialized site ordering called ``seq(BF)'' \cite{ye2026practical}; in 3D, it differs from a conventional QTT chain and more effectively captures the dominant correlations associated with finite-difference derivatives and the directional coupling in Maxwell's equations. This makes CTN a more suitable representation for the anisotropic operators and boundary effects that arise in FDTD with PML.

In this paper, we propose a CTN-augmented FDTD algorithm for efficient 2D and 3D EM simulations with PML. The proposed CTN-FDTD method represents the permittivity, electric and magnetic fields, and the associated auxiliary PML variables using CTNs. The finite-difference operators are applied directly in the tensor-compressed form, without forming the full 2D or 3D arrays throughout the time-stepping process. A central ingredient of this approach is the use of rank-reduction strategies to maintain compactness over long simulations. In particular, local smoothing and global smoothing are applied to suppress rank inflation caused by repeated derivative updates and noise accumulation. These steps preserve the dominant physical content while removing spurious short-scale components that would otherwise increase the tensor ranks\rev{; they also change the solution to some extent, in particular the smoothing of the permittivity changes the material, so we quantify the smoothing errors separately from the truncation errors}. The resulting CTN-FDTD algorithm retains the essential accuracy of standard FDTD while substantially reducing memory and CPU costs.

The proposed method is assessed on a range of 2D and 3D benchmarks involving wave propagation, current excitation, and heterogeneous material distributions. \rev{In 3D, the scattering by a dielectric sphere is compared with standard FDTD and, at matched accuracy, with QTT chains that use the same implementation.} The results show that CTN-FDTD substantially reduces memory and CPU complexity relative to standard FDTD while preserving the accuracy and stability of the underlying solution. In the 3D setting, the method is particularly compelling: CTN-FDTD can handle a large grid with 35 trillion spatial grid points, which would require 281 TB storage for one 3D array with standard FDTD. These examples highlight the practical potential of CTN compression for high-resolution EM modeling and establish the method as a scalable alternative to standard and QTT-based FDTD algorithms.

The rest of the paper is organized as follows. Section~\ref{sec:fdtd} reviews the standard 2D split-field FDTD and 3D ADE-FDTD formulations, which motivate the tensor-compressed approach. Section~\ref{sec:cnt_fdtd} introduces the CTN-FDTD algorithm, including the CTN representation, the smoothing strategies, and the complexity analysis. Section~\ref{sec:numerical} presents the numerical experiments and performance validation for representative large-scale electromagnetic simulations.

	\section{The Standard FDTD Algorithms with PML}\label{sec:fdtd}	
In this work, we are interested in simulating Maxwell equations in a lossless, unbounded domain with constant permeability: 
\begin{eqnarray}
\frac{\partial \mathbf{H}}{\partial t}
= -\frac{1}{\mu_0}\nabla \times \mathbf{E},
\quad
\label{eq:dHdt}
\frac{\partial \mathbf{E}}{\partial t}
= \frac{1}{\varepsilon}\nabla \times \mathbf{H}
- \frac{1}{\varepsilon}\mathbf{J}.
\label{eq:dEdt}
\end{eqnarray}
where $\textbf{E}$ denotes the electric fields, $\textbf{H}$ denotes the magnetic fields, $\textbf{J}$ denotes the electric current density, $\mu_0$ is the free-space permeability, and $\varepsilon$ is the (not necessarily constant) permittivity of the domain. To illustrate the applicability of QTT to various types of FDTD algorithms and dimensionalities \cite{atef2009finite}, we briefly summarize the 2D FDTD algorithm with split-field PML and 3D FDTD algorithm with auxiliary differential equation (ADE)-PML below. \rev{We chose these two commonly used PML formulations for illustration; the CTN-FDTD algorithm applies to other PML formulations as well, such as the uniaxial PML \cite{wang2025ahfdtd} or the complex-frequency-shifted PML \cite{hou2025cdhie}. These formulations also update the fields and their auxiliary variables by Hadamard products with coefficients that vary only across the PML layers, shifts and sums, so the choice of the PML formulation does not change the complexity of CTN-FDTD and should not notably affect its cost.} For simplicity, we assume the computational domain (including the PML region) is discretized into a $N^d$ grid (with $N_{pml}$ PML cells per boundary) with $d=2,3$ and spatial resolution $dx=dy=dz=ds$.    
\subsection{2D FDTD with split-field PML}
For 2D $\textrm{TM}_z$ scattering, the standard FDTD algorithm with split-field PML \cite{atef2009finite} splits the $E_{z}$ field as $E_{z}=E_{zx}+E_{zy}$. The update equations at time steps $n+\frac{1}{2}$ and $n+1$ for each $N\times N$ field quantity are
\begin{eqnarray}
H_x^{\,n+\frac{1}{2}}
&=
G\odot H_x^{\,n-\frac{1}{2}}
-
H\odot\rev{\delta^{+}_y}\!\left(E_{zx}^{\,n}+E_{zy}^{\,n}\right),\label{eq:2dHx}
\\
H_y^{\,n+\frac{1}{2}}
&=
A\odot H_y^{\,n-\frac{1}{2}}
+
B\odot \rev{\delta^{+}_x}\!\left(E_{zx}^{\,n}+E_{zy}^{\,n}\right),
\\
E_{zx}^{\,n+1}
&=
C\odot E_{zx}^{\,n}
-
D\odot \rev{\delta^{-}_y} H_x^{\,n+\frac{1}{2}}
-
\frac{dt}{2\varepsilon_0}\,J_z^{\,n+\frac{1}{2}},\label{eq:2dJz1}
\\
E_{zy}^{\,n+1}
&=
E \odot E_{zy}^{\,n}
+
F\odot \rev{\delta^{-}_x} H_y^{\,n+\frac{1}{2}}
-
\frac{dt}{2\varepsilon_0}\,J_z^{\,n+\frac{1}{2}}.\label{eq:2dJz2}
\end{eqnarray}
Here \rev{$[\delta^{+}_x u]_{ij} = (u_{i+1,j} - u_{ij})/ds$ and $[\delta^{-}_x u]_{ij} = (u_{ij} - u_{i-1,j})/ds$, $i,j=0,\ldots,N-1$, are the forward and backward finite-difference approximations of the $1^{st}$ order spatial derivative along $x$ ($\delta^{\pm}_y$ are defined likewise along $y$); differences that would reach outside the grid are set to zero (they lie at the outer edge of the PML).} $\odot$ denotes Hadamard product. The coefficient matrices are: 
\begin{align}
&G_{ij} = \frac{\mu_0 - \tfrac{1}{2}dt\,\sigma_{x,ij}^{*}}{\mu_0 + \tfrac{1}{2}dt\,\sigma_{x,ij}^{*}}, \, H_{ij} = \frac{\rev{dt}}{\mu_0 + \tfrac{1}{2}dt\,\sigma_{x,ij}^{*}},
\nonumber\\
&A_{ij} = \frac{\mu_0 - \tfrac{1}{2}dt\,\sigma_{y,ij}^{*}}{\mu_0 + \tfrac{1}{2}dt\,\sigma_{y,ij}^{*}}, \, B_{ij} = \frac{\rev{dt}}{\mu_0 + \tfrac{1}{2}dt\,\sigma_{y,ij}^{*}},
\nonumber\\
&C_{ij} = \frac{\varepsilon_{ij} - \tfrac{1}{2}dt\,\sigma_{x,ij}}{\varepsilon_{ij} + \tfrac{1}{2}dt\,\sigma_{x,ij}}, \, D_{ij} = \frac{\rev{dt}}{\varepsilon_{ij} + \tfrac{1}{2}dt\,\sigma_{x,ij}},
\nonumber\\
&E_{ij} = \frac{\varepsilon_{ij} - \tfrac{1}{2}dt\,\sigma_{y,ij}}{\varepsilon_{ij} + \tfrac{1}{2}dt\,\sigma_{y,ij}}, \, F_{ij} = \frac{\rev{dt}}{\varepsilon_{ij} + \tfrac{1}{2}dt\,\sigma_{y,ij}},\label{eq:A-H}
\end{align}	
where $\sigma_x$ and $\sigma_y$ represent the fictitious conductivity profiles in the PML along $x$ and $y$ directions, respectively. The starred sigma are multiplied by the squared wave impedance, $\sigma_{c}^{*} = \frac{\mu_0}{\varepsilon}\,\sigma_{c}$, $c=x,y$. 
% $\sigma_x^{*} = \frac{\mu_0}{\varepsilon}\,\sigma_x$, and $\sigma_y^{*} = \frac{\mu_0}{\varepsilon}\,\sigma_y$. 
Note that for grid points outside the PML, $N_{pml}\leq i,j<N-N_{pml}$, $A_{ij}=G_{ij}=C_{ij}=E_{ij}=1$, \rev{$B_{ij}=H_{ij}=dt/\mu_0$ and $D_{ij}=F_{ij}=dt/\varepsilon_{ij}$.} In \eqref{eq:2dJz1} and \eqref{eq:2dJz2}, it's assumed that the current density resides in free-space. Consider the number of time steps to be $N_t=O(N)$, the computational cost and memory requirement for the standard 2D FDTD scale as $O(N^3)$ and $O(N^2)$, respectively.

\subsection{3D FDTD with ADE-PML}
For fully 3D electromagnetics modeling, the ADE-PML algorithm \cite{atef2009finite} introduces a set of 12 auxiliary field quantities $\Psi_{h,xy}$, $\Psi_{h,xz}$, $\Psi_{h,yx}$, $\Psi_{h,yz}$, $\Psi_{h,zx}$, $\Psi_{h,zy}$, $\Psi_{e,xy}$, $\Psi_{e,xz}$, $\Psi_{e,yx}$, $\Psi_{e,yz}$, $\Psi_{e,zx}$, and $\Psi_{e,zy}$, which together with the 6 electric and magnetic field quantities $E_x,E_y,E_z,H_x,H_y,H_z$, are discretized into 18 $N\times N\times N$ mode-3 tensors. Note that simpler versions of the PML requiring fewer auxiliary fields can also be considered \cite{jiang2018simple}. For notational simplicity, define $\mathbf{P}^{\,n+\frac{1}{2}}_{h}
=[-\Psi_{h,xy}^{\,n+\frac{1}{2}}+\Psi_{h,xz}^{\,n+\frac{1}{2}};\Psi_{h,yx}^{\,n+\frac{1}{2}}-\Psi_{h,yz}^{\,n+\frac{1}{2}};-\Psi_{h,zx}^{\,n+\frac{1}{2}}+\Psi_{h,zy}^{\,n+\frac{1}{2}}]$ and $\mathbf{P}^{\,n+1}_{e}=\left[\Psi_{e,xy}^{\,n+1}-\Psi_{e,xz}^{\,n+1};-\Psi_{e,yx}^{\,n+1}+\Psi_{e,yz}^{\,n+1};\Psi_{e,zx}^{\,n+1}-\Psi_{e,zy}^{\,n+1}\right]$. The update equations at time steps $n+\frac{1}{2}$ and $n+1$ are as follows:

\begin{align}
\Psi_{h,ac}^{\,n+\frac{1}{2}}
&=
g_{h,c}\odot\Psi_{h,ac}^{\,n-\frac{1}{2}}
+
f_{h,c}\odot\rev{\delta^{-}_c} E_{b}^{\,n},\label{eq:3dpsih}\\
\mathbf{H}^{\,n+\frac{1}{2}}
&=
\mathbf{H}^{\,n-\frac{1}{2}}
-
dt\,\mu_0^{-1}
\left(
\rev{\nabla^{-}_h} \times \mathbf{E}^{\,n}
-
\mathbf{P}_h^{\,n+\frac{1}{2}}\right),\\
\Psi_{e,ac}^{\,n+1}
&=
g_{e,c}\odot\Psi_{e,ac}^{\,n}
+
f_{e,c}\odot\rev{\delta^{+}_c} H_{b}^{\,n+\frac{1}{2}},\\
\mathbf{E}^{\,n+1}
&=
\mathbf{E}^{\,n}
+
\frac{dt}{\varepsilon_0 \varepsilon_r}
\odot\left(
\rev{\nabla^{+}_h} \times \mathbf{H}^{\,n+\frac{1}{2}}
+
\mathbf{P}_e^{\,n+1}\right)-
\frac{dt}{\varepsilon_0}\,\mathbf{J}^{\,n+\frac{1}{2}},\label{eq:3dEupdate}
\end{align}
where the subscripts $(a,b,c)=(y,z,x)$, $(z,y,x)$, $(x,z,y)$, $(z,x,y)$, $(x,y,z)$ or $(y,x,z)$. Here $\varepsilon_0$ is the free-space permittivity, $\varepsilon_r=\varepsilon/\varepsilon_0$ is the relative permittivity, and 
\begin{eqnarray}
\left[g_{p,c} \right]_{ijk} =& \rev{\frac{\varepsilon_0}{\varepsilon_0 + dt \left( \left[\sigma_{pc}\right]_{ijk} + \left[\alpha_{pc}\right]_{ijk} \right)}}\\
\left[f_{p,c}\right]_{ijk} =& \rev{\frac{-dt \, \left[\sigma_{pc}\right]_{ijk}}{\varepsilon_0 + dt \left( \left[\sigma_{pc}\right]_{ijk} + \left[\alpha_{pc}\right]_{ijk} \right)}} \, .\label{eq:ab}	
\end{eqnarray}
Subscript $p$ denotes parameters corresponding to either the electric fields (`$e$') or the magnetic fields (`$h$').
% \begin{eqnarray}
% a_{ec,ijk} =& \frac{-dt \sigma_{ec,ijk}}{ds\left( \varepsilon_0 + dt \left( \sigma_{ec,ijk} + \alpha_{ec,ijk} \right) \right)}\\
% b_{ec,ijk} =& \frac{\varepsilon_0}{\varepsilon_0 + dt \left( \sigma_{ec,ijk} + \alpha_{ec,ijk} \right)}\label{eq:ab}	
% \end{eqnarray}
$\sigma_{pc}$ and $\alpha_{pc}$ are the conductivity and shift parameters in the PML regions along $c=x,y,z$ direction \rev{(for $p=h$, the magnetic conductivity and shift scaled by $\varepsilon_0/\mu_0$, i.e., the electric profiles sampled at the magnetic-field locations)} \cite{atef2009finite}. Note that \rev{$\left[\sigma_{pc}\right]_{ijk}=\left[\alpha_{pc}\right]_{ijk}=0$ outside the PML regions, where $f_{p,c}=0$: the auxiliary fields remain zero there and \eqref{eq:3dpsih}--\eqref{eq:3dEupdate} reduce to the standard Yee updates.} 
%Similarly $b_{mc}$ and $a_{mc}$ can be computed from the conductivity $\sigma_{mc}$ and shift $\alpha_{mc}$ used for the H update. 
In \eqref{eq:3dpsih} to \eqref{eq:3dEupdate}, \rev{$\nabla^{\pm}_h$ represent the discretized curl operators, specifically $\nabla^{\pm}_h \times \mathbf{X}=[\delta^{\pm}_y X_z - \delta^{\pm}_z X_y;\delta^{\pm}_z X_x - \delta^{\pm}_x X_z;\delta^{\pm}_x X_y - \delta^{\pm}_y X_x]$}  	
where \rev{the forward and backward differences $[\delta^{+}_x u]_{ijk} = (u_{i+1,j,k} - u_{ijk})/ds$ and $[\delta^{-}_x u]_{ijk} = (u_{ijk} - u_{i-1,j,k})/ds$ (and likewise along $y$ and $z$) are finite-difference approximations of the $1^{st}$ order spatial derivative.} Just like \eqref{eq:2dJz1}, here we assume $\mathbf{J}$ only nonzero in the background media. Moreover, one can add impressed E and H fields, $\mathbf{E}^{\,n+1}_{imp}$ and $\mathbf{H}^{\,n+1}_{imp}$, to $\mathbf{E}^{\,n+1}$ and $\mathbf{H}^{\,n+1}$, instead of using $\mathbf{J}^{n+\frac{1}{2}}$ in \eqref{eq:3dEupdate}. 

Consider the number of time steps to be $N_t=O(N)$, the computational cost and memory requirement for the standard 3D FDTD scale as $O(N^4)$ and $O(N^3)$, respectively. It's worth mentioning that in addition to the 18 $N\times N\times N$ tensors for $\Psi_{e,ac}$, $\Psi_{h,ac}$, $\mathbf{E}$ and $\mathbf{H}$, this algorithm also stores 12 $N\times N\times N$ tensors for the coefficients $f_{e,c}$, $g_{e,c}$, $f_{h,c}$ and $g_{h,c}$, and one $N\times N\times N$ tensor for $1/\varepsilon_r$, leading to a total of at least 31 mode-3 tensors. This memory requirement becomes prohibitively high for large $N$. For example, $N=2048$ leads to a memory requirement of at least 2.13 TB, which can be significantly reduced by tensor compression.    
	
\section{CTN-FDTD Algorithms with PML}\label{sec:cnt_fdtd}	
In this section, we describe the proposed CTN augmentation of the 2D and 3D FDTD algorithms in \cref{sec:fdtd}, which effectively brings down the computation and memory requirement to $O(R^3 N\log N)$ and $O(R^2 \log N)$\rev{ in 2D and to $O(R^4N+R^3N\log N)$ and $O(R^3+R^2\log N)$ in 3D (\Cref{sec:cc})}, where $R$ is the CTN rank. Assume $N=2^L$ in both 2D and 3D simulations. We will describe the notations and operations mostly in terms of the 3D FDTD algorithm, and briefly mention their adaptation to 2D FDTD algorithms. 
\subsection{CTN representation of FDTD quantities}\label{sec:cnt}		
\begin{figure}[!t]
	\centering
	\includegraphics[width=0.8\columnwidth]{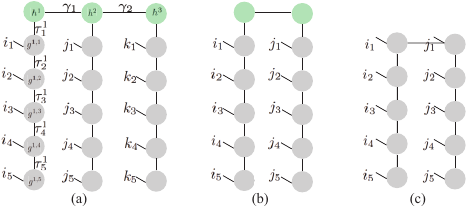}
	\caption{Tensor diagrams for (a) comb tensor network for a $N\times N\times N$ tensor with $N=2^5$, (b) comb tensor network for a $N\times N$ tensor with $N=2^5$, and (c) QTT for a $N\times N$ tensor with $N=2^5$ \label{fig:qtt}}
\end{figure}

For any $N\times N\times N$ tensor $T$ encountered in the 3D FDTD algorithm, the comb tensor network \cite{Dolgov2013comb, chepiga2019comb,ye2024quantized}, also closely related to the extended TT or QTT-Tucker decomposition \cite{Dolgov2013comb}, starts with a TT representation as 
\begin{align}
T_{i, j, k} &= \sum_{\gamma_{1}, \gamma_{2}}s^1_{i,\gamma_{1}}s^2_{j,\gamma_{1},\gamma_{2}}s^3_{k,\gamma_{2}},
	\label{eq:TT_rep_simple}
\end{align}
where $i,j,k=0,\ldots,N-1$ denote the grid points along each spatial dimension, and $\gamma_1=0,\ldots, R_1-1$ and $\gamma_2=0,\ldots, R_2-1$ are bond indices with bond dimensions (i.e., numerical ranks) $R_1$ and $R_2$. Next, each of the cores $ s^1_{i,\gamma_{1}}s^2_{j,\gamma_{1},\gamma_{2}}s^3_{k,\gamma_{2}}$ with fixed $\gamma_{1},\gamma_{2}$ is reshaped into a mode-$L$ tensor via 
   	\begin{equation}
	\label{eq:multi-index}
	i=\sum_{\ell=1}^{L}2^{L-\ell}i_{\ell}, \,
    j=\sum_{\ell=1}^{L}2^{L-\ell}j_{\ell}, \,
	k=\sum_{\ell=1}^{L}2^{L-\ell}k_{\ell},
\end{equation}
with $i_\ell,j_\ell,k_\ell\in\{0,1\}$. They are represented as QTTs, 
\begin{eqnarray}
	s^1_{i_1, \dots, i_L, \gamma_1} = \sum_{\tau^1_1, \dots, \tau^1_{L}}h^1_{\gamma_1,\tau^1_1}\Big(\prod_{\ell = 1}^{L
	}g^{1,\ell}_{i_\ell, \tau^1_{\ell},\tau^1_{\ell+1}} \Big)\label{eq:QTT_rep_simple_S1}\\
	s^2_{j_1, \dots, j_L, \gamma_1,\gamma_2} = \sum_{\tau^2_1, \dots, \tau^2_{L}}h^2_{\gamma_1,\gamma_2,\tau^2_1}\Big(\prod_{\ell = 1}^{L
	}g^{2,\ell}_{j_\ell, \tau^2_{\ell},\tau^2_{\ell+1}} \Big)\label{eq:QTT_rep_simple_S2}\\
	s^3_{k_1, \dots, k_L, \gamma_2} = \sum_{\tau^3_1, \dots, \tau^3_{L}}h^3_{\gamma_2,\tau^3_1}\Big(\prod_{\ell = 1}^{L
}g^{3,\ell}_{k_\ell, \tau^3_{\ell},\tau^3_{\ell+1}} \Big).	
	\label{eq:QTT_rep_simple_S3}
\end{eqnarray}
Here, $\tau^p_{\ell}$, $p=1,2,3$, $\ell=1,\ldots,L$ are the bond indices with bond dimensions $R^p_{\ell}$, and $\tau^p_{L+1}$ are dummy bond indices of size 1. Putting together \eqref{eq:TT_rep_simple}, \eqref{eq:QTT_rep_simple_S1}, \eqref{eq:QTT_rep_simple_S2} and \eqref{eq:QTT_rep_simple_S3}, the comb tensor network representation of $T$ is 
\begin{align}
	&
    T_{i_1, \dots, i_L, j_1, \dots, j_L, k_1, \dots, k_L} =
    \nonumber
    \\ 
    & \hspace{1.5cm}
    \sum_{\gamma_{1}, \gamma_{2}}s^1_{i_1, \dots, i_L,\gamma_{1}}s^2_{j_1, \dots, j_L,\gamma_{1},\gamma_{2}}s^3_{k_1, \dots, k_L,\gamma_{2}}.
	\label{eq:comb_rep_simple}
\end{align} 
We refer to $h^p$ (subscripts dropped) as the spine tensor cores, and $g^{p,\ell}$ (subscripts dropped) as the branch tensor cores. For example, \Cref{fig:qtt}(a) shows the tensor diagram of a comb tensor network for a $N\times N\times N$ tensor with $N=2^5$. The 3 spine tensor cores are shown in green and the 15 branch tensor cores are shown in gray. It's worth noting that for mode-3 tensors (in 3D FDTD), the CTN representation is structurally different compared with QTT. In contrast, a 2-dimensional CTN (see \Cref{fig:qtt}(b)) is essentially equivalent to the QTT representation (see \Cref{fig:qtt}(c)) of % Specifically, the CTN is equivalent to the TT representation in \eqref{eq:TT_rep_simple} for 
the reshaped tensor $T_{i_L,\ldots,i_1,j_1,\ldots,j_L}$. 
Note that the indices of $T$ are permuted such that $i_L,\ldots,i_1$ represents the least significant bit to most significant bit order, and $j_1,\ldots,j_L$ represents the most significant bit to least significant bit order. This is referred to as the ``seq(BF)" ordering for QTT in \cite{ye2026practical}. 

\begin{figure*}[!t]
	\centering
	\includegraphics[width=\columnwidth]{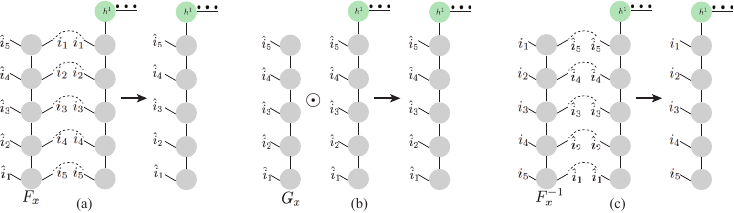}
	\caption{Tensor diagrams for the three steps of global smoothing of a CTN $T$ into $S$. (a) Forward DFT: $\mathrm{vec}(\hat{T})=\big(F_z\otimes F_y \otimes F_x\big)\mathrm{vec}(T)$. (b) Application of the filter: $\mathrm{vec}(\hat{S})=\big(G_z\otimes G_y \otimes G_x\big)\odot\mathrm{vec}(\hat{T})$. (c) Inverse DFT: $\mathrm{vec}(S)=\big(F_z^{-1}\otimes F_y^{-1} \otimes F_x^{-1}\big)\mathrm{vec}(\hat{S})$. Only the action on the first branch is shown, and the spine cores remain unchanged. \label{fig:filter}}
\end{figure*}	

\subsection{CTN operations in FDTD}\label{sec:cnt_op}		
Assume that all 31 field quantities in 3D FDTD are represented in the CTN format (see \Cref{sec:cnt_init} regarding details of their initialization). The update of the 12 field quantities in \eqref{eq:3dpsih} to \eqref{eq:3dEupdate}, i.e., $\Psi_{e,ac}$, $\Psi_{h,ac}$, $\mathbf{E}$ and $\mathbf{H}$,  requires the following CTN operations: 
\subsubsection{Hadamard product}
The Hadamard (element-wise) product of two CTNs of same dimensions, $\tilde{T}=T\odot \hat{T}$, can be expressed as~\cite{Dolgov2013comb} 
\begin{eqnarray}
\tilde{h}^p=h^p\otimes \hat{h}^p,~\tilde{g}^{p,\ell}_{q,:,:}={g}^{p,\ell}_{q,:,:}\otimes \hat{g}^{p,\ell}_{q,:,:}
\end{eqnarray}
where $p=1,2,3$, $q=i_\ell,j_\ell,k_\ell$, $\ell=1,\ldots,L$ and $\otimes$ is the Kronecker product of two matrices or tensors. The exact Hadamard product of two CTNs leads to increased ranks $\tilde{R}_1=R_1\hat{R}_1$, $\tilde{R}_2=R_2\hat{R}_2$, and $\tilde{R}_l^p=R_l^p\hat{R}_l^p$, and a SVD-based rank truncation of the exact result is oftentimes applied to reduce its CTN ranks~\cite{Dolgov2013comb}. The cost of the rank truncation procedure can be further reduced by using the zip-up algorithm on each of the TT branches \cite{paeckel2019time,stoudenmire2010minimally}.

\subsubsection{Summation}
The summation of two CTNs of same dimensions, $\tilde{T}=T+ \hat{T}$, can be performed as~\cite{Dolgov2013comb}:
\begin{eqnarray}
	\tilde{h}^p=h^p\oplus \hat{h}^p,~\tilde{g}^{p,\ell}_{q,:,:}={g}^{p,\ell}_{q,:,:}\oplus \hat{g}^{p,\ell}_{q,:,:},~\tilde{g}^{p,L}_{q,:,1}=[{g}^{p,L}_{q,:,1}, \hat{g}^{p,L}_{q,:,1}]
\end{eqnarray}
where $p=1,2,3$, $q=i_\ell,j_\ell,k_\ell$, $\ell=1,\ldots,L-1$. The exact summation of two CTNs leads to increased ranks $\tilde{R}_1=R_1+\hat{R}_1$, $\tilde{R}_2=R_2+\hat{R}_2$, and $\tilde{R}_l^p=R_l^p+\hat{R}_l^p$, and a SVD-based rank truncation of the exact result is applied to reduce its ranks. Also, since the multiplication of a CTN with a scalar is trivial, the operation $\tilde{T}=\alpha T+ \beta\hat{T}$ can be performed in a similar fashion. 

\subsubsection{Kronecker product of QTTs}\label{sec:kron}
Let $T^p$ with $p=1,2,3$ denote a tensorized vector along the specified spatial dimension. It can be expressed in QTT form, $T^p_i=\prod_{\ell = 1}^{L
}g^{p,\ell}_{i^{\ell}, \tau^p_{\ell},\tau^p_{\ell+1}}$, where $\tau^p_{1}$ and $\tau^p_{L+1}$ are dummy indices of size one. The CTN representation of their Kronecker product, $T=T^1\otimes T^2\otimes T^3$, uses the QTT cores as the branch cores, and the spine cores are set to 1.
% $h^1_{:,:}=h^2_{:,:,:}=h^3_{:,:}=1$, $R_1=R_2=1$.  
\subsubsection{Shift}
The finite-difference-based spatial derivatives, \rev{$\delta^{\pm}_x$, $\delta^{\pm}_y$ and $\delta^{\pm}_z$}, require shifting the mode-3 tensor along one specific mode. For example, \rev{$\delta^{\pm}_x$} can be rewritten as 
 \begin{eqnarray}
 \rev{\delta^{+}_x = \frac{1}{ds}(S^+\otimes I \otimes I - I \otimes I \otimes I ),}\nonumber\\
 \rev{\delta^{-}_x = \frac{1}{ds}(I \otimes I \otimes I - S^-\otimes I \otimes I )}\label{eq:shift}
 \end{eqnarray}
 where the shift \rev{operators are $(S^{+}f)_i=f_{i+1}$ and $(S^{-}f)_i=f_{i-1}$} where $i=0,\ldots,N-1$ for a length-$N$ vector $f$\rev{, with $f_N=f_{N-1}$ and $f_{-1}=f_0$ so that the differences vanish at the ends of the grid}. \rev{$S^{\pm}$} can be efficiently represented as \rev{QTTs} as described in~\cite{kazeev_multilevel_2013}. When applied to a CTN, the \rev{shift terms} of \eqref{eq:shift} \rev{amount} to updating only the branch cores in \eqref{eq:QTT_rep_simple_S1} with the QTT representation of \rev{$S^{\pm}$}.

\rev{\subsubsection{Rounding of the updates}\label{sec:rounding}
Each update in \eqref{eq:3dpsih}--\eqref{eq:3dEupdate} is a short sum of Hadamard products and shifted fields; e.g., the update of $H_x$ sums $H_x^{\,n-\frac{1}{2}}$, the shifted copies of $E_y$ and $E_z$ and two auxiliary fields. Instead of forming such a sum term by term with an SVD-based truncation after every addition, the implementation used for the 3D scattering examples in \Cref{sec:sphere,sec:cube} truncates the whole sum once by a randomized rounding on the comb tree \cite{daas2023randomized}: the terms are contracted with Gaussian random sketches of twice the maximum rank (1.75 times on the spine bonds), the sketched network is orthogonalized from the leaves to the spine, and all bonds are then truncated together, to the tolerance and the maximum rank, from the Gram matrices of the orthogonalized network \cite{grasedyck2010hierarchical}.}

\subsection{CTN Initialization in FDTD}\label{sec:cnt_init}		
In order to maintain low memory consumption across the entire 3D FDTD simulation, the time-independent auxiliary arrays (\rev{$f_{e,c}$, $g_{e,c}$, $f_{h,c}$ and $g_{h,c}$}; and $1/\varepsilon_r$ in \eqref{eq:3dEupdate}) and impressed source and field quantities ($\mathbf{J}^{n+\frac{1}{2}}$ in \eqref{eq:3dEupdate}, $\mathbf{E}^{\,n+1}_{imp}$ and $\mathbf{H}^{\,n+1}_{imp}$) must be constructed as CTNs without ever fully forming any $N\times N\times N$ tensors in 3D FDTD. The same principle applies to the $N\times N$ tensors in 2D FDTD with $A,B,C,D,E,F,G$ of \eqref{eq:A-H} and $J_z^{n+\frac{1}{2}}$ of \eqref{eq:2dJz2}. 

The auxiliary quantities are often separable with small tensor ranks. To see this, consider \rev{$f_{e,x}$} of \eqref{eq:ab} in 3D FDTD for an example. We have \rev{$\left[f_{e,x}\right]_{ijk}=T^x_iT^y_jT^z_k$} where $T^y_j=T^z_k=1$, and $T^x_{N_{pml}:N-N_{pml}-1}=0$ since $\left[\sigma_{ex}\right]_{ijk}=0$ for $i=N_{pml}, \ldots, N-N_{pml}-1$. The QTT representation of $T^x$ has ranks of at most $2N_{pml}$. The CTN representation of  \rev{$f_{e,x}=T^x\otimes I\otimes I$} (described in \Cref{sec:kron}) requires at most $O(N)$ computation and peak storage. The same result holds for other auxiliary quantities \rev{$f_{e,c}$, $g_{e,c}$, $f_{h,c}$ and $g_{h,c}$}. 

For the dielectric constant, 
we only consider profiles such that its inverse can be represented as a sum of a few separable components, $1/\varepsilon_r=\sum_{i} T^x\otimes T^y\otimes T^z$ . Thus its CTN representation requires $O(N)$ computation and peak storage. For more complex permittivity profiles, one can leverage the TT-cross-type algorithms \cite{oseledets2010tt} to approximate $1/\varepsilon_r$ as a CTN with low complexity. 

For the excitation, we only consider current sources $\mathbf{J}^{n+\frac{1}{2}}$ (or impressed fields $\mathbf{E}^{\,n+1}_{imp}$ and $\mathbf{H}^{\,n+1}_{imp}$) in the form of $\mathbf{J}^{n+\frac{1}{2}}=f_nT^x\otimes T^y\otimes T^z$. Here $f_n$ is the discrete time signature. $T^c$ with $c=x,y,z$ is a standard basis vector, an all-one vector, or a localized Gaussian vector. Therefore, the CTN representation of the excitation function can be constructed with $O(\log N)$ computation and peak storage costs.  

\subsection{Smoothing Techniques}\label{sec:smooth}
For smooth, separable, and/or localized tensors, such as the auxiliary arrays, current source, and impressed field quantities, their CTN ranks are small constant. However, the updates in \eqref{eq:3dpsih} to \eqref{eq:3dEupdate} can lead to significant rank increases in $\Psi_{h,ac}^{\,n+\frac{1}{2}}$, $\Psi_{e,ac}^{\,n+1}$, $\mathbf{E}^{\,n}$ and $\mathbf{H}^{\,n+\frac{1}{2}}$ due to both accumulation of numerical errors and complex field phenomena involving inhomogeneous $\epsilon_r$\rev{; the former was analyzed for TT-FDTD in \cite{zhou2025tensor}, when the fields become small after their absorption in the PML}. We propose two types of computationally efficient smoothing techniques for the field quantities and $\epsilon_r$, which represent practical strategies to trade off simulation fidelity with computational efficiency. \rev{Both change the solution of the discretized problem. The local smoothing adds artificial dissipation, and the global smoothing of $1/\epsilon_r$ changes the material: a sharp interface becomes a graded one, which removes part of the physical reflection and diffraction at the interface. These smoothing errors are separate from the truncation error and do not decrease with the truncation tolerance. We therefore quantify them separately in \Cref{sec:numerical}: against standard FDTD with and without the same local smoothing (\Cref{tab:wg}), and for scatterers with sharp interfaces, whose permittivity is not smoothed (\Cref{sec:sphere}).}

\subsubsection{Local smoothing}\label{sec:local_smooth}   
For any field quantity $T_n$ in the left hand side (LHS) of \eqref{eq:3dpsih} to \eqref{eq:3dEupdate} in 3D FDTD and \eqref{eq:2dHx} to \eqref{eq:2dJz2} in 2D FDTD, we perform a local smoothing operation after their update:  
\begin{align}
	T_n \leftarrow T_n \!+\! (-1)^{\frac{m}{2}+1}\eta\rev{\sum_{c}\left(D_c\right)^{\frac{m}{2}}T_n,}\label{eq:local_smooth}
\end{align}
\rev{where $D_c$ is the centered second difference along $c=x,y,z$, $[D_x u]_{ijk}=u_{i+1,j,k}-2u_{ijk}+u_{i-1,j,k}$ (and likewise along $y$ and $z$). Hence $(D_c)^{m/2}=ds^m\,\partial_c^m+O(ds^{m+2})$ is the centered $m^{th}$ difference, with the stencils $(1,-2,1)$ for $m=2$ and $(1,-6,15,-20,15,-6,1)$ for $m=6$. Near the two ends of each axis, which lie in the PML, the stencil needs values outside the grid; they are extrapolated by the polynomial of degree $m+1$ through the $m+2$ values nearest to the boundary\rev{ in the 2D examples and the 3D example without scatterer. In the 3D scattering examples (\Cref{sec:sphere,sec:cube}), they are instead mirror images of the values inside the grid. This makes $(D_c)^{m/2}$ symmetric with zero column sums, so the smoothing conserves the mean of each field (the extrapolated closure does not, which slowly injected a spatially constant field in our 3D tests)}.} Equation \eqref{eq:local_smooth} amounts to adding artificial dissipation terms into \eqref{eq:dHdt}. \rev{The dissipation coefficient $\eta\geq0$ is dimensionless and controls the strength of the dissipation. For example, at PPW=19.2 (the 3D examples) the $6^{th}$-order smoothing with $\eta=10^{-3}$ attenuates the propagating wave by $1.2\times10^{-6}$ per time step (about $0.5\%$ over 4000 steps), while it damps grid-scale components by up to $19\%$ per step.} $m$ is an even integer, and in this paper we consider $m=2$ (referred to as $2^{nd}$ order local smoothing) and $m=6$ (referred to as  $6^{th}$ order local smoothing) \cite{ye2026practical}. 

We note that these smoothing operations can be performed efficiently given that $T_n$ is in the CTN format, as they involve only scalar product, summation, and shifting operations on CTN, whose implementation has been explained in \Cref{sec:cnt_op}.  

\subsubsection{Global smoothing}\label{sec:global_smooth}   
For the same field quantities as those in \Cref{sec:local_smooth}, one can also use an alternative global smoothing technique. Also, we apply this global smoothing operation to $1/\varepsilon_r$ \cite{nguyen2026tensor,nguyen2026tensorFDTD}\rev{, which turns sharp interfaces into graded ones and thus changes the scatterer (see \Cref{sec:smooth}); among the examples in \Cref{sec:numerical}, only the cube in \Cref{sec:cube} uses it}. Unlike existing global smoothing operations \cite{nguyen2026tensor}\rev{\cite{nguyen2026ttfdtd3d}} that suffer from $O(N^d)$ computational and storage costs, we propose an efficient global smoothing operation in this paper. \rev{It acts directly on the CTN of the quantity to be smoothed, regardless of how that CTN was obtained (from a sum of separable terms, as in our examples, or by TT-cross \cite{oseledets2010tt} for general profiles), so the full $N^d$ array is never formed. In contrast, the Gaussian interface smoothing of \cite{nguyen2026ttfdtd3d} applies a 3D FFT to the material array on the full grid before the compression. The QTT construction of the DFT is due to \cite{chen2025direct}; what we add is its application to the branches of a CTN.} Specifically, consider a mode-3 tensor $T$ to be smoothed. The proposed global smoothing operation is defined as follows:
\begin{eqnarray}
\mathrm{vec}(\hat{T})=\big(F_z\otimes F_y \otimes F_x\big)\mathrm{vec}(T),\label{eq:fft}\\
\mathrm{vec}(\hat{S})=\big(G_z\otimes G_y \otimes G_x\big)\odot\mathrm{vec}(\hat{T}),\label{eq:filter}\\
\mathrm{vec}(S)=\big(F_z^{-1}\otimes F_y^{-1} \otimes F_x^{-1}\big)\mathrm{vec}(\hat{S}),\label{eq:ifft}
\end{eqnarray} 
where $F_x=F_y=F_z$ is the discrete Fourier transform (DFT) matrix with $F_{x,\hat{i}i}=\exp(-2\pi\mathrm{j}i\hat{i}/N)$, $i,\hat{i}=0,\ldots,N-1$. $G_x=G_y=G_z$ is the Gaussian low-pass filter with $G_{x,\hat{i}}=\exp(-{2|\hat{i}|^2}/{(N^2k_c^2)})$ if $0\leq\hat{i}<N/2$ and $G_{x,\hat{i}}=\exp(-{2|\hat{i}-N|^2}/{(N^2k_c^2)})$ if $N/2\leq\hat{i}<N$, and $k_c$ is the normalized cutoff frequency. Alternatively, we can consider the Sigmoid low-pass filter with $G_{x,\hat{i}}=
[1+\exp((2|\hat{i}|/N-k_c)/t)]^{-1}$ if $0\leq\hat{i}<N/2$ and $G_{x,\hat{i}}=
[1+\exp((2|\hat{i}-N|/N-k_c)/t)]^{-1}$ if $N/2\leq\hat{i}<N$, where $k_c$ is the normalized cutoff frequency and $t=0.1$ controls the transition width of the filter.

The smoothed quantity only needs to computed once before being saved into memory. However, when working at high resolutions, it is desirable to perform the calculation efficiently without fully forming any $N\times N \times N$ tensor. 
In what follows, we describe the global smoothing operation of \eqref{eq:fft} to \eqref{eq:ifft}. Let $\hat{i}_\ell=0,1$ and $\ell=1,\ldots,L$ be the multi-indices for $\hat{i}$ as in \eqref{eq:multi-index}. First, the Gaussian filter $G_x,G_y,G_z$ can be represented as a QTT with at most $O(N)$ memory and computational costs, before the time marching starts. Next, the 1D DFT matrix $F_x,F_y,F_z$ can be directly represented as a QTT, in $O(\log N)$ computation and storage, with explicitly defined tensor cores \cite{chen2025direct}. Specifically, one can write 
\begin{align}
	F_{x,\hat{i}_1,\ldots,\hat{i}_L,i_1,\ldots,i_L}
	&=
	\sum_{\tau_2,\ldots,\tau_L}
	\prod_{\ell=1}^{L}
	f^{\ell}_{i_\ell,\hat{i}_{N-\ell+1},\tau_\ell,\tau_{\ell+1}},\nonumber\label{eq:qtt_FFT}
	\\
	f^{1}_{a,b,0,d}
	&=
	\sum_{c=0}^{R_F-1} k(a,b,c,d),\nonumber
	\\
	f^{\ell}_{a,b,c,d}
	&=
	k(a,b,c,d),
	\qquad 1<\ell\leq L,\nonumber
	\\
	k(a,b,c,d)
	=&
	P^c\left(\frac{c+\alpha_d}{2}\right)
	\exp\left(-\pi\mathrm{j}b(a+\alpha_d)\right).
\end{align}
Here, $\tau_1=\tau_L=0$, $\tau_\ell=0,\ldots,R_F-1$, where $R_F$ is the QTT rank. $a,b=0,1$ and $c,d=0,\ldots,R_F-1$. $P^c$ are degree-$R_F$ Lagrange interpolating polynomials on $[0,1]$, and $\alpha_{d} = (1 - \cos(\pi d/R_F))/{2}$ are the Chebyshev-Lobatto nodes. Note that for each QTT core, $i_\ell$ is paired with $\hat{i}_{N-\ell+1}$, i.e., the output uses the bit-reversal ordering to ensure small constant QTT ranks \cite{chen2025direct}. As the QTT cores are explicitly constructed, they require the optimal $O(\log N)$ computation and storage units. Moreover, the inverse DFT $F_x^{-1}$ can be expressed as a QTT in a similar way to \eqref{eq:qtt_FFT}.  

Next, the three steps of global smoothing in the CTN format are illustrated in \Cref{fig:filter} with $L=5$. In the forward DFT step of \eqref{eq:fft}, the QTT $F_c$, $c=x,y,z$ is independently applied to each branch of the CTN via contraction along the indices $i_\ell$, $\ell=1,\ldots,5$ without modifying the spine cores. In the filter application step of \eqref{eq:filter}, the QTT $G_c$, $c=x,y,z$ is independently applied to each branch of the CTN via the Hadamard product without modifying the spine cores. The inverse DFT step of \eqref{eq:ifft} is similar to the forward DFT step. After the three steps, an additional SVD-based rank truncation of $S$ is performed to keep the CTN ranks low.      

\subsection{Complexity Summary}\label{sec:cc}
In this subsection, we analyze the computational and memory complexities of CTN-FDTD for 2D and 3D simulations. We assume that FDTD uses $N^d$ grid points and $N_t=O(N)$ time steps, and the CTN ranks for time-independent auxiliary quantities, current sources and impressed fields are small constant. Also, we use 
$R=\max_n\{\max\{R_1,R_2,R_\ell^p\}\}$, $p=1,2,3$, $\ell=1,\ldots,L$, $n=0,\ldots,N_t-1$, 
hence $R$ to denote the maximum rank of all bonds in the CTN across all time steps and all fields. 
%for $\Psi_{h,ac}^{\,n+\frac{1}{2}}$, $\Psi_{e,ac}^{\,n+1}$, $\mathbf{E}^{\,n}$, $\mathbf{H}^{\,n+\frac{1}{2}}$ in 3D-FDTD, and $H_x^{\,n+\frac{1}{2}}$, $H_y^{\,n+\frac{1}{2}}$, $E_{zx}^{\,n+1}$, $E_{zy}^{\,n+1}$ in 2D-FDTD. 
As we will see in \Cref{sec:numerical}, $R$ is not necessarily constant independent of $N$ and needs to be accounted for in the complexity analysis. \rev{The bounds below are therefore conditional on the compressibility of the fields, not a general complexity guarantee for electromagnetic scattering. $R$ also depends on $N$ differently when the discretization of a fixed physical problem is refined (\Cref{fig:scaling_point_Jz}(a)--(b)) and when the electrical size grows at a fixed PPW (\Cref{fig:scaling_point_Jz}(c)--(d) and \Cref{fig:cc_3d_planewave}).}
		
First, the initialization of time-independent auxiliary quantities and time dependent excitation, as explained in \Cref{sec:cnt_init}, require $O(dN)$ and $O(d\log N)$ computational costs, respectively. The $O(N)$ costs come from the temporary formation of $O(N)$-sized arrays\rev{, also for the Gaussian global smoothing of $1/\varepsilon_r$. A nonseparable coefficient, such as the sphere in \Cref{sec:sphere}, is compressed from its full array at $O(N^d)$ cost. TT-cross \cite{oseledets2010tt} could reduce these costs to $O(d\log N)$}.

Next, the update operations in \Cref{sec:cnt_op} and the smoothing operations in \Cref{sec:smooth} are fully performed in the CTN format\rev{. With bond dimensions of at most $R$, the $dL$ branch cores have at most $2R^2$ entries each; the spine cores are $R\times R$ matrices in 2D, but in 3D the middle spine core $h^2$ in \eqref{eq:QTT_rep_simple_S2} has three bonds and up to $R^3$ entries. A field quantity therefore takes $O(R^d+dR^2\log N)$ storage.

Each update and each smoothing combines field quantities with operators whose ranks do not depend on $N$: sums of at most five terms, Hadamard products with coefficients of small CTN rank, shifts, and the smoothing operators of \Cref{sec:smooth}. Formed exactly, the result has bond dimensions of at most $\kappa R$, where $\kappa$ does not depend on $N$. Its SVD-based truncation costs $O(d\kappa^3R^3\log N)$ for the QR factorizations and SVDs of the $dL$ branch cores, which are $2\kappa R\times\kappa R$ matrices, and, in 3D, $O(\kappa^4R^4)$ for the middle spine core, reshaped into a $(\kappa R)^2\times\kappa R$ matrix. A time step therefore costs $O(R^{d+1}+dR^3\log N)$. The randomized rounding of \Cref{sec:rounding} avoids the bonds of dimension $\kappa R$ but has the same order, as it contracts the middle spine core at $O(R^4)$ cost per term in 3D. With $N_t=O(N)$ time steps, the total cost is $O(R^{d+1}N+dR^3N\log N)$, which is $O(N\log N)$ only as long as $R$ stays bounded.} These costs are in stark contrast to the $O(N^d)$ per time step computational cost \rev{($O(N^{d+1})$ in total)} and the $O(N^d)$ storage cost of standard FDTD, as shown in \Cref{tab:cc}.

It is worth mentioning that some of our numerical experiments visualize the $\textbf{E}$ and $\textbf{H}$ fields at certain time steps. In 2D FDTD the $N\times N$ tensors are fully reconstructed from the CTN format, and in 3D FDTD only $N\times N$ slices are reconstructed from the CTN format. This requires an additional $O(N^2)$ computational and memory costs\rev{ per visualized time step}, as listed in \Cref{tab:cc}. However, we treat this as the optional postprocessing step and do not count this towards the overall complexity of the CTN-FDTD algorithms. \rev{The DFT of $\mathbf{E}$ on a plane for the analysis in \Cref{sec:sphere} likewise adds $O(N^2)$ work per time step.}

 \begin{table} %[!tp]
	\centering	
	\setlength\tabcolsep{3pt}
	\footnotesize
	\begin{tabular}{|c|c|c|c|}
		\hline
		\multirow{2}{*}{Operations} & \multicolumn{2}{c|}{\rev{Computation}} & \multirow{2}{*}{Memory} \\
		\cline{2-3}
		& \rev{per time step} & \rev{total} & \\
		\hline
		Initialization & \rev{--} & \rev{$O(dN)$, $O(N^d)$$^*$} & \rev{$O(dN)$, $O(N^d)$$^*$} \\
		\hline
		Update & \rev{$O(C_R)$} & \rev{$O(NC_R)$} & \rev{$O(M_R)$} \\
		\hline
		Smoothing & \rev{$O(C_R)$} & \rev{$O(NC_R)$} & \rev{$O(M_R)$} \\
		\hline
		Visualization & \rev{--} & $O(N^{2})$ & $O(N^{2})$ \\
		\hline
		Standard FDTD & \rev{$O(N^{d})$} & $O(N^{d+1})$ & $O(N^{d})$ \\
\hline
	\end{tabular}
	\caption{CPU and memory complexities of CTN-FDTD with $N^d$ grid points and $N_t=O(N)$ time steps. $R$ represents the maximum CTN ranks of the field quantities during time marching, which is not necessarily constant\rev{; the logarithmic per-step costs assume that $R$ stays bounded. $C_R=R^{d+1}+dR^3\log N$ and $M_R=R^d+dR^2\log N$, i.e., $O(R^3\log N)$ and $O(R^2\log N)$ in 2D, and $O(R^4+R^3\log N)$ and $O(R^3+R^2\log N)$ in 3D}. Visualization refers to reconstruction of the full $N\times N$ field quantities in 2D or $N\times N$ slices in 3D for the plotting purpose. \rev{$^*$: $O(N^d)$ only for the sphere in \Cref{sec:sphere}.} but can be removed by TT-cross.}\label{tab:cc}
\end{table}

	\section{Numerical Results}\label{sec:numerical}
	This section presents several 2D and 3D examples to demonstrate the accuracy, efficiency and scalability of the proposed CTN-FDTD algorithm. Both the proposed CTN-FDTD and the standard FDTD algorithms are implemented in Python and executed on one CPU node of the Perlmutter machine at NERSC in Berkeley. The node has two $64$-core AMD EPYC 7763 processors and a total of 512 GB of 2133MHz DDR4 memory. 
	
	\subsection{CTN-FDTD for 2D EM Simulations}
	In what follows, we investigate the scalability, accuracy and applicability of CTN-FDTD for various excitations in homogeneous and inhomogeneous 2D domains.  
	\subsubsection{Point Current Source}
	\rev{Our canonical test is a challenging point current source excitation.} The computational domain is $128\times 128 ~\mu m^2$ with a current source $J_z(\mathbf{r},t)=\sin(2\pi ft)\delta(\mathbf{r}-\mathbf{r}_0)$ with $\mathbf{r}_0$ being the domain center. As we will see, the radiation pattern of a point source with high frequencies is not low rank \cite{ye2026practical}. Therefore, smoothing techniques are crucial in order to reduce the computational and memory costs. In what follows, we first investigate the effects of smoothing, then compare CTN with QTT in terms of rank behaviors, and finally show the scaling plots.    
	
\begin{figure}[!t]
	\centering
	\includegraphics[width=\columnwidth]{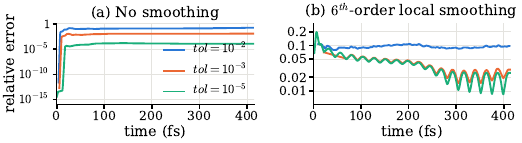}
	\caption{	
		Relative errors vs. time between CTN-FDTD (of different compression tolerances) and standard FDTD for a 2D point source with $L=7$ ($N=128$), $N_t=175$ and PPW=14.2. (a) No smoothing. \rev{(b)} $6^{th}$-order local smoothing with $\eta=0.01$. \label{fig:error_point_Jz}}
\end{figure}

\begin{figure}[!t]
		\centering
		\includegraphics[width=\columnwidth]{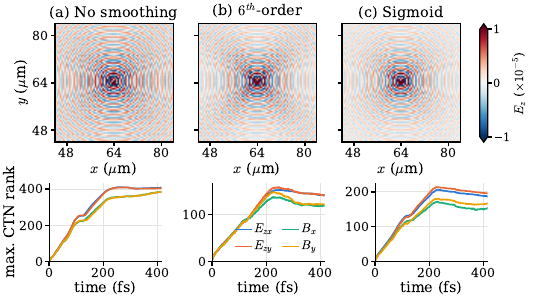}
		\caption{	
			The field pattern \rev{$E_z$ in the central $40\times40~\mu m^2$ of the domain} at $t=213 fs$ \rev{(top)} and the \rev{maximum CTN rank of each field quantity} over time \rev{(bottom)} for a 2D point source with $L=9$\rev{ and} $f=169.7$ THz with (a) no smoothing, \rev{(b)} $6^{th}$-order local smoothing ($\eta=0.01$) and \rev{(c)} Sigmoid global smoothing ($k_c=1.2$). \label{fig:fix_ppw_wavelength_point_Jz}}
	\end{figure}

	First, we consider $f=21.213$ THz, $L=7$ ($N=128$, $ds=1~\mu m$), $dt=ds/(\sqrt{2}c)$, and $N_t=175$. This corresponds to PPW=14.2. The relative errors at each time step, compared with standard FDTD ($L=7$), with no smoothing and $6^{th}$-order local smoothing ($\eta=0.01$) are shown in \Cref{fig:error_point_Jz}. Although no smoothing (\Cref{fig:error_point_Jz}(a)) can achieve the desired compression error $10^{-2}, 10^{-3}, 10^{-5}$, the use of $6^{th}$-order local smoothing can achieve at most $10^{-2}$ accuracy. Therefore, we use at most $tol=10^{-3}$ in the following tests. \rev{The reason is that the error against standard FDTD then consists of two parts: the truncation error, which is about $10\,tol$ without smoothing (\Cref{fig:error_point_Jz}(a)), and the change of the solution by the smoothing itself, which does not depend on $tol$ and is about 1--9\% here after the first 25~fs (\Cref{fig:error_point_Jz}(b)). For $tol\leq10^{-3}$, the smoothing error dominates, so the curves for $tol=10^{-3}$ and $10^{-5}$ nearly coincide, and only $tol=10^{-2}$ adds a visible truncation error.

This suggests how we select the parameters. The tolerance (or the maximum rank) only needs to make the truncation error smaller than the smoothing error, since a smaller tolerance raises the ranks, and the cost grows with the ranks as in \Cref{sec:cc}, without improving the accuracy. The smoothing coefficient $\eta$ and the order $m$ trade the smoothing error against the ranks: a wave resolved by PPW points per wavelength is attenuated by about $\eta(2\pi/\mathrm{PPW})^m$ per time step, so a higher order and a smaller $\eta$ change the solution less, while the stronger damping of grid-scale components keeps the ranks lower. In practice, we therefore choose $\eta$, within the stability bound $\eta\leq4/(3\,d\,4^{m/2})$, as small as the ranks and, in 3D, the stability of the simulations at the maximum rank allow. The photonic crystal waveguide below quantifies this balance: $\eta=10^{-2}$ and $10^{-3}$ give errors of 33--37\% and 5.5--5.9\% against standard FDTD, with 2--4 times less storage than standard FDTD and twice as much, respectively (\Cref{tab:wg}).} For $L=9$ \rev{and} $f=169.7$ THz (PPW=\rev{7.07}), the field patterns at time step $700$ and CTN ranks at each time step, subject to no smoothing, $6^{th}$-order local smoothing ($\eta=0.01$) and Sigmoid global smoothing ($k_c=1.2$), are shown in \Cref{fig:fix_ppw_wavelength_point_Jz}. On one hand, the CTN ranks approach their theoretical maximum of $2^L=512$ if no smoothing is used, but are significantly reduced by various smoothing techniques \rev{(to at most about 160--215)}. On the other hand, the use of smoothing can lead to distortion of the field patterns to some extent (as illustrated by \Cref{fig:error_point_Jz} as well). Overall, $6^{th}$-order local smoothing seems to yield the best tradeoff between accuracy and low-rankness. It is worth mentioning that the ranks saturate after the wave front reaches the PML boundaries at around $t=213~fs$, as seen by all \rev{three} tests.        
	
\begin{figure}[!t]
	\centering
	\includegraphics[width=\columnwidth]{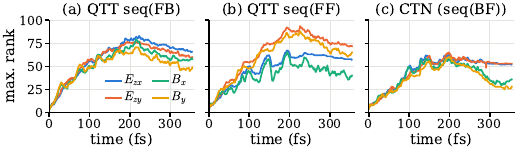}
	\caption{\rev{The maximum ranks of the field quantities} over time for the point current source with $L=10$, $f=21.213$ THz, $tol=10^{-3}$, and Sigmoid global smoothing. (a) QTT with seq(FB) ordering, (b) QTT with seq(FF) ordering, and (c) QTT with seq(BF) ordering, i.e., CTN. \label{fig:ordering_point_Jz}}
\end{figure}

	Next, we compare CTN with QTT in terms of rank behaviors using $L=10$, $f=21.213$ THz, $tol=10^{-3}$, and Sigmoid global smoothing. We test QTT with different orderings, specifically ``seq(FB)", ``seq(FF)" and ``seq(BF)" orderings in \cite{ye2026practical}. Note that the CTN is equivalent to ``seq(BF)" ordering for 2D, but distinct from QTT in 3D (see \Cref{fig:qtt}). The ranks over time are plotted in \Cref{fig:ordering_point_Jz} where one can see that CTN gives the lowest ranks. In all other experiments of this paper, we only use CTN instead of QTT\rev{, except for the 3D comparison with QTT chains in \Cref{sec:sphere,sec:cube}}.

\begin{figure}[!t]
	\centering
	\includegraphics[width=\columnwidth]{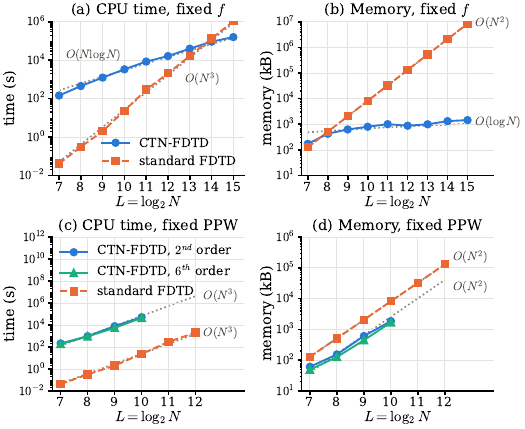}
	\caption{Scaling results for the 2D point current source using CTN-FDTD and standard FDTD. Top: (a) The CPU time and (b) memory usage when $f=21.213$ THz for $L=7$ to $15$. Bottom: (c) The CPU time and (d) memory usage when PPW=\rev{7.07} for $L=7$ to $12$. \label{fig:scaling_point_Jz}}
\end{figure}

	Finally, we investigate the CPU and memory scaling of CTN-FDTD for the point current source test, by fixing $f$ and PPW, respectively. When fixing $f=21.213$ THz \rev{(a refinement of a fixed physical problem)}, we use a $tol=10^{-3}$ and $2^{nd}$-order local smoothing. The smoothing coefficients are chosen as $\eta=0.016(2^{7-L})$. The total run time and memory requirement of CTN-FDTD and standard FDTD for $L=7$ to $15$ are shown in \rev{\Cref{fig:scaling_point_Jz}}(a)-(b). The maximum CTN ranks for \rev{$L=7$} to 15 are 24, 35, 46, 53, 70, 78, 90, 109, 121, and CTN-FDTD exhibits the desired $O(N\log N)$ CPU and $O(\log N)$ memory scaling. CTN-FDTD outperforms standard FDTD in terms of memory and CPU when $L>8$ and $L>13$, respectively. 
	
	When fixing PPW=\rev{7.07} \rev{(an increase of the electrical size at a fixed PPW)}, we use $tol=10^{-3}$ with $2^{nd}$-order and $6^{th}$-order local smoothing ($\eta=0.01$). Even with smoothing, it is expected that $R=O(N)$ in \Cref{tab:cc} in the high-frequency regime and CTN-FDTD loses advantages over standard FDTD. This can be observed from \rev{\Cref{fig:scaling_point_Jz}}(c) and (d). Although there is still a constant ratio of memory reduction compared with standard FDTD, the CPU time is not comparable with standard FDTD.

	\subsubsection{Photonic Crystal Waveguide }		
\begin{figure*}[!t]
	\centering
	\includegraphics[width=0.99\textwidth]{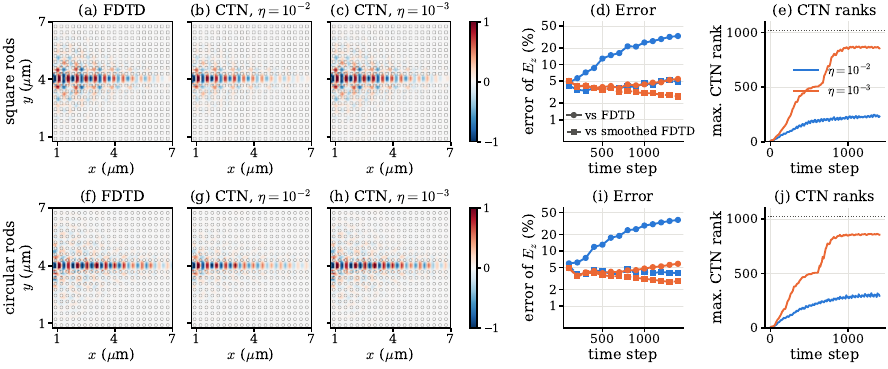}
	\caption{\rev{Results for the 2D photonic crystal waveguide ($L=10$, $N=1024$, PPW=16.64, $N_{pml}=100$, rods not smoothed). Top: dielectric rods with a square cross section. Bottom: dielectric rods with a circular cross section. (a)--(c), (f)--(h): $E_z$ at time step 1400 from standard FDTD and from CTN-FDTD with the $2^{nd}$-order local smoothing, $\eta=10^{-2}$ and $10^{-3}$, normalized by the 99.5th percentile of $|E_z|$ of standard FDTD (rod outlines in gray). (d), (i): relative error of $E_z$ outside the PML over time, against standard FDTD (solid) and against standard FDTD with the same smoothing (dashed). (e), (j): maximum CTN rank of the four field quantities over time; the dotted line is the maximum possible rank, 1024.}\label{fig:wg}}
\end{figure*}

Finally we leverage CTN-FDTD to simulate a realistic 2D photonic crystal structure as shown in \Cref{fig:wg}. The computational domain is $7.68\times 7.68 ~\mu m^2$ with a $20\times 20$ array of electric rods of $\epsilon_r=11.56$. Each rod's cross section is either a circle of radius $76.875~nm$ or a square of side length $153.75~nm$. The current source is  $J_z(x,y,t)=\delta(x-x_0)W(y)
\cos\left(\frac{\pi(y-y_0)}{a}\right)\sin(2\pi f t)$, where $x_0=0.7575~\mu m$, $y_0=3.9786~\mu m$, $a=0.3075~\mu m$, $W(y)=1$ when $|y-y_0|<a/2$ and 0 otherwise. The frequency is $f=707.1$ THz. The domain is discretized as $L=10$ ($N=1024$), $ds=7.5~nm$ (PPW=16.64), $dt=0.01768~fs$, $N_{pml}=100$, and $N_t=1400$. 

\rev{CTN-FDTD uses a rank truncation tolerance $tol=10^{-3}$ and the $2^{nd}$-order local smoothing of the fields with $\eta=10^{-2}$ or $10^{-3}$. The permittivity is not smoothed: the grid points inside the rods have $\epsilon_r=11.56$, as in standard FDTD. We compare each CTN-FDTD simulation with standard FDTD and, to separate the error of the tensor compression from that of the smoothing, also with standard FDTD that applies the same local smoothing \eqref{eq:local_smooth} after every update. \Cref{fig:wg} shows $E_z$ at time step 1400, the relative errors of $E_z$ outside the PML over time, and the CTN ranks over time, and \Cref{tab:wg} lists the ranks, storage, run times and errors at time step 1400.

The two smoothing strengths trade accuracy for compression. The compression error alone, against standard FDTD with the same smoothing, is 2.7--4.8\% for both strengths and both rod shapes, and it changes little after the first 100 time steps (\Cref{fig:wg}(d) and (i)). With $\eta=10^{-2}$, the maximum CTN ranks are 249 (square rods) and 320 (circular rods), and the four field quantities take 8.8 and 15.8~MB at time step 1400, against 33.6~MB for standard FDTD. However, the smoothing itself changes the solution by 33--37\%: inside the rods, a wavelength spans only 4.9 cells, and the smoothing with $\eta=10^{-2}$ attenuates such a wave by 1.4\% per time step. With $\eta=10^{-3}$, CTN-FDTD is within 5.5\% (square rods) and 5.9\% (circular rods) of standard FDTD, but the CTN ranks grow to about 870, close to the maximum of 1024, the field quantities take 70~MB, twice the storage of standard FDTD, and the simulations take 6.0 (square rods) and 3.4 (circular rods) times as long as with $\eta=10^{-2}$. At the same tolerance, the fields of standard FDTD at time step 1400 have ranks of only 89--126 (square rods), so the larger CTN ranks mostly carry accumulated truncation errors, which the stronger smoothing removes together with part of the solution. The circular rods cost more than the square rods because their $1/\epsilon_r$ is not separable in $x$ and $y$ (CTN rank 58 against 28).}

\begin{table}[!t]
	\centering
	\setlength\tabcolsep{3pt}
	\footnotesize
	\begin{tabular}{|c|c|c|c|c|c|c|}
		\hline
		\multirow{2}{*}{Rods} & \multirow{2}{*}{$\eta$} & Max. & Storage & \multirow{2}{*}{Time} & \multicolumn{2}{c|}{Error of $E_z$ (\%) vs.} \\
		\cline{6-7}
		& & rank & (MB) & & FDTD+sm. & FDTD \\
		\hline
		Square & $10^{-2}$ & 249 & 8.8 & 1.33~h & 4.8 & 33.0\\
		Square & $10^{-3}$ & 873 & 69.5 & 7.97~h & 2.7 & 5.5\\
		\hline
		Circular & $10^{-2}$ & 320 & 15.8 & 8.59~h & 4.0 & 36.7\\
		Circular & $10^{-3}$ & 865 & 70.0 & 29.4~h & 2.8 & 5.9\\
		\hline
		\multicolumn{2}{|c|}{Standard FDTD} & -- & 33.6 & 65~s & \multicolumn{2}{c|}{--}\\
		\hline
	\end{tabular}
	\caption{\rev{CTN-FDTD for the photonic crystal waveguide ($L=10$, $tol=10^{-3}$, 1400 time steps): maximum CTN rank over time, storage of the four field quantities at time step 1400, run time, and error of $E_z$ outside the PML at time step 1400 against standard FDTD with the same local smoothing (FDTD+sm.) and without it (FDTD).}}\label{tab:wg}
\end{table}

	\subsection{CTN-FDTD for 3D EM Simulations}
	In this section, we present \rev{three} 3D examples of homogeneous and inhomogeneous domains subject to plane-wave excitation to demonstrate the efficacy and scalability of CTN-FDTD for 3D EM simulations\rev{: a plane wave without scatterer on a grid of up to $32768^3$ points, a dielectric sphere on a grid of $512^3$ points, whose CTN-FDTD solution is compared with standard FDTD and with QTT chains, and a dielectric cube on a grid of $4096^3$ points}.
	\subsubsection{Plane-wave Excitation without Scatterer}\label{sec:pw3d}	
	In this example, we consider a plane wave propagating in the homogeneous media. The computational domain is $1\times 1\times 1 ~k m^3$ with an impressed plane wave field $
	\mathbf{E}^{\mathrm{imp}}(\mathbf{r},t)=\hat{\mathbf{e}}\,
	S(t)\sin(2\pi ft)\delta(\hat{\mathbf{k}}\cdot\mathbf{r}-s_0)$ with $\hat{\mathbf{k}}=(1,0,0)$, $\hat{\mathbf{e}}=(0,1,0)$, $s_0=(N_{pml}+1)ds$, $f$ is the excitation frequency, and $S(t)$ is a smooth envelope such that $S(t)=0$ for $t>9/f$.  
	
\begin{figure}[!t]
	\centering
	\includegraphics[width=\columnwidth]{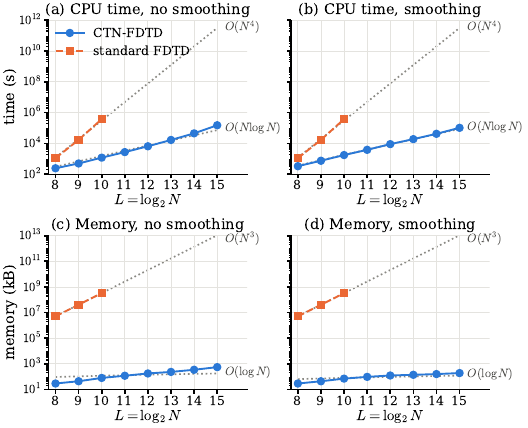}
	\caption{Scaling results for the 3D plane wave excitation without scatterer using CTN-FDTD and standard FDTD for $L=8$ to $15$ ($N^3=16.7$ million to $35$ trillion), and PPW=19.2. Left: (a) The CPU time and (c) memory usage when CTN-FDTD applies no smoothing. Right: (b) The CPU time and (d) memory usage when CTN-FDTD applies $6^{th}$-order local smoothing with $\eta=0.001$.\label{fig:cc_3d_planewave}}
\end{figure}

	\rev{We} fix PPW=19.2 and vary $f$ from 4 MHz to 512 MHz\rev{, i.e., we increase the electrical size at a fixed PPW}. This corresponds to $L=8$ to $15$, and $N_t=500$ to $64000$. We tested CTN-FDTD \rev{with $tol=10^{-4}$,} with both $\eta=0.001$ and no smoothing. It is remarkable that $L=15$ corresponds to $N^3=32768^3=35$ trillion grid points. It would require 281 TB memory to store a single field quantity by the standard FDTD algorithm. As a result, we can only run standard FDTD up to $L=10$, $N^3\approx 1$ billion. The total run time and memory requirement of CTN-FDTD and standard FDTD are shown in \Cref{fig:cc_3d_planewave}. It is worth noting that without smoothing, the CTN ranks will increase slightly  for plane wave field patterns, but adding the smoothing keeps the CTN ranks constant as $L$ increases. Overall, at $L=15$, CTN-FDTD can achieve a projected 10 orders of magnitude reduction in memory and 6 orders of magnitude speedup in computation time, compared with standard FDTD.

\subsubsection{\rev{Dielectric Sphere: Comparison with Standard FDTD and QTT}}\label{sec:sphere}
\begin{figure*}[!t]
	\centering
	\includegraphics[width=0.99\textwidth]{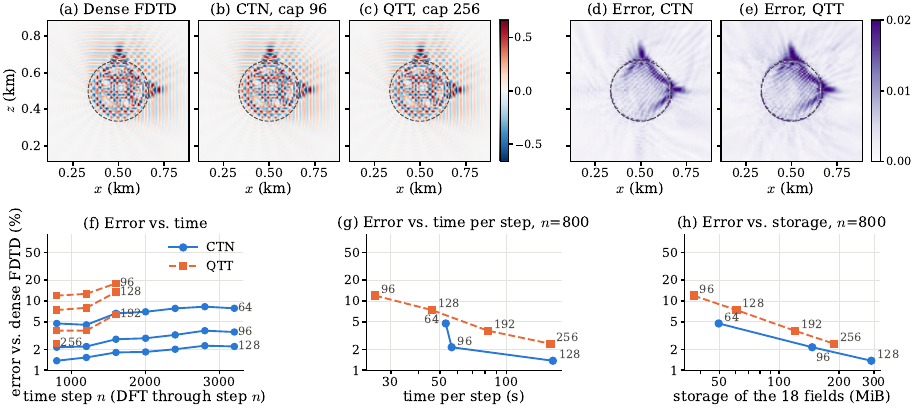}
	\caption{\rev{3D plane-wave excitation with a dielectric sphere for $f=8$~MHz and $L=9$ (PPW=19.2): CTN-FDTD and QTT chains against standard FDTD. (a)--(c) $\mathrm{Re}(E_{s,y})$ on the plane $y=0.5~km$ (DFT at $f$ through time step $n=800$, normalized by the maximum of (a)) from standard FDTD, CTN-FDTD with the maximum rank 96, and QTT with the maximum rank 256; the dashed circle is the sphere. (d)--(e) $|\mathbf{E}_s-\mathbf{E}_s^{\mathrm{FDTD}}|$ of the CTN and QTT solutions in (b)--(c), relative to the maximum of $|\mathbf{E}_s^{\mathrm{FDTD}}|$. (f) Error against standard FDTD vs. $n$, labeled by the maximum rank. (g) Error vs. time per step and (h) error vs. storage of the 18 field quantities for $n=800$.}}
	\label{fig:sphere}
\end{figure*}

\rev{We consider the same $1\times 1\times 1 ~k m^3$ domain with a dielectric sphere of radius $166~m$ and $\epsilon_r=2$ at its center, illuminated by two plane waves as in \Cref{sec:pw3d} with $f=8$~MHz and $(\hat{\mathbf{k}},\hat{\mathbf{e}})=((1,0,0),(0,1,1)/\sqrt{2})$ and $((0,0,1),(1,1,0)/\sqrt{2})$. With $L=9$ ($N=512$, $ds=1.953~m$, $dt=3.676~ns$, PPW=19.2, $N_{pml}=56$), the radius is 85 cells ($ka=27.8$); $N_t=3200$. The permittivity is not smoothed (a staircase interface), and its CTN is compressed from the full $512^3$ array (ranks 167 on the branches and 84 on the spine, error $4\times10^{-12}$). CTN-FDTD uses $tol=10^{-5}$, maximum ranks of 64, 96 and 128, and the $6^{th}$-order local smoothing of $\mathbf{E}$ and $\mathbf{H}$ with $\eta=3\times10^{-3}$. The reference is standard FDTD with the identical discretization and no compression or smoothing (20 minutes, 33~GB). The error is the relative 2-norm, on the plane $y=0.5~km$ outside the PML, of the difference of the scattered fields $\mathbf{E}_s$: the running DFTs $\sum_{n'\le n}\mathbf{E}^{n'}e^{2\pi\mathrm{j}fn'dt}\,dt$, minus those without the sphere.

With the maximum ranks 64, 96 and 128, which take 53.2, 56.5 and 161.9~s per time step, the errors are 4.7\%, 2.2\% and 1.4\% for $n=800$, when the pulses have just passed the sphere, and 7.8\%, 3.6\% and 2.2\% at the end (\Cref{tab:sphere}, \Cref{fig:sphere}(a), (b), (d) and (f)); they are largest near the surface where the waves exit. The errors include the local smoothing, which attenuates the waves by at most 0.3\% outside and 2.2\% inside the sphere through step 800: the truncation dominates for the ranks 64 and 96, while for 128 the smoothing may account for a comparable part. Standard FDTD itself differs from the Mie series by 33\% at this resolution due to numerical dispersion and staircase approximation of the sphere.

We also run the same code with every quantity represented as a QTT chain over all $3L=27$ bits and all other settings unchanged, in the sequential ordering $x_1,\ldots,x_L,y_1,\ldots,y_L,z_1,\ldots,z_L$ ($x_1$ the most significant bit; the QTT of $1/\epsilon_r$ has ranks up to 743 and error $5.6\times10^{-4}$). Since equal maximum ranks do not give equal accuracy, we raise the QTT rank until the errors match (\Cref{tab:sphere}, \Cref{fig:sphere}(c) and (e)--(h)). With the rank 256, the QTT (2.4\%) is still less accurate than CTN with the rank 96 (2.2\%), at 2.8 times its time per step and 1.3 times its storage (about 3 times the time at equal accuracy). The QTT errors also grow faster in time. Only the spine bonds and the first two bonds of each branch of the CTN reach the maximum rank, against 9--14 of the 26 QTT bonds.}

\begin{table}[!t]
	\centering
	\setlength\tabcolsep{3pt}
	\footnotesize
	\begin{tabular}{|c|c|c|c|c|c|c|}
		\hline
		\multirow{2}{*}{Network} & Max. & Storage & Time/step & \multicolumn{3}{c|}{Error (\%), $n=$} \\
		\cline{5-7}
		& rank & (MiB) & (s) & 800 & 1200 & 1600 \\
		\hline
		CTN & 64 & 49.5 & 53.2 & 4.74 & 4.52 & 6.62\\
		CTN & 96 & 146.5 & 56.5 & 2.15 & 2.20 & 2.80\\
		CTN & 128 & 290.9 & 161.9 & 1.37 & 1.53 & 1.81\\
		\hline
		QTT & 96 & 37.0 & 25.5 & 11.9 & 12.7 & 18.0\\
		QTT & 128 & 60.8 & 46.0 & 7.39 & 7.96 & 13.4\\
		QTT & 192 & 119.6 & 82.4 & 3.72 & 3.73 & 6.42\\
		QTT & 256 & 188.5 & 157.8 & 2.41 & -- & --\\
		\hline
	\end{tabular}
	\caption{\rev{CTN and QTT chains for the dielectric sphere ($L=9$, $tol=10^{-5}$, one node): storage of the 18 field quantities, time per step (mean over steps 801--1000, all fields at the maximum rank), and error of the scattered field against standard FDTD for the DFT through step $n$. ``--'': not run.}}\label{tab:sphere}
\end{table}

\subsubsection{Plane-wave Excitation with \rev{a Cubic} Scatterer}\label{sec:cube}
\begin{figure*}[!t]
	\centering
	\includegraphics[width=0.99\textwidth]{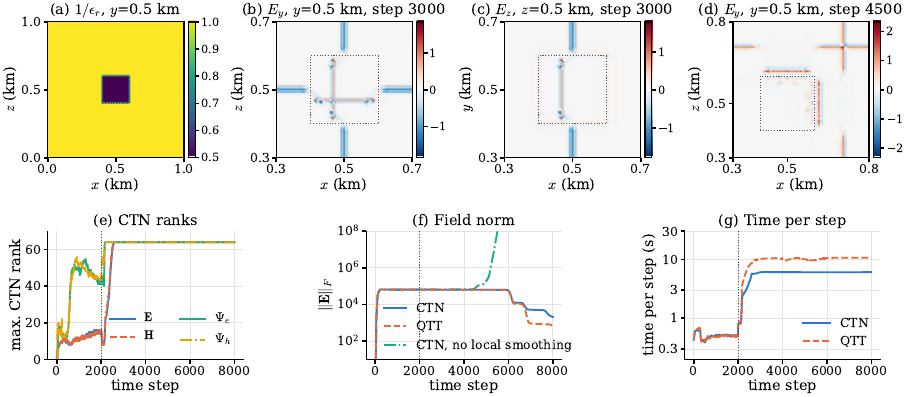} 
	\caption{\rev{3D plane-wave excitation with a cube scatterer using CTN-FDTD with the maximum rank 64 for $f=64$~MHz and $L=12$ ($N=4096$, PPW=19.2). (a) Smoothed $1/\epsilon_r$ on the plane $y=0.5~km$. (b) $E_y$ on the plane $y=0.5~km$ and (c) $E_z$ on the plane $z=0.5~km$ at time step 3000, and (d) $E_y$ on the plane $y=0.5~km$ at time step 4500, near the cube (dotted). (e) Maximum CTN ranks of $\mathbf{E}$, $\mathbf{H}$, $\Psi_{e,ac}$ and $\Psi_{h,ac}$ per time step. (f) Norm of $\mathbf{E}$ for CTN with and without the $6^{th}$-order local smoothing ($\eta=10^{-3}$ from time step 2000, dotted line) and for QTT with the same settings. (g) Time per step of CTN and QTT (running median over 51 steps, visualization steps excluded).}}
	\label{fig:L12_3d_planewave_cube}
\end{figure*}

\rev{We consider the same domain with a cube of side $200~m$ and $\epsilon_r=2$ at its center, illuminated by the plane waves of \Cref{sec:sphere} with $f=64$~MHz, with $L=12$ ($N^3=68.7$ billion grid points, $ds=0.2441~m$, $dt=0.4596~ns$, PPW=19.2, $N_{pml}=448$) and $N_t=8000$; standard FDTD would need 17~TB for its 31 arrays. $1/\epsilon_r$ is smoothed by the Gaussian global smoothing with $k_c=0.02$, which grades the interface over a standard deviation of 16 cells (0.83 wavelength), so the smoothed cube is a different, graded scatterer (CTN ranks of $1/\epsilon_r$ at most 5); sharp interfaces are assessed with the sphere in \Cref{sec:sphere}. CTN-FDTD uses $tol=10^{-5}$, the maximum rank 64, and the $6^{th}$-order local smoothing of $\mathbf{E}$, $\mathbf{H}$ and $\Psi_{e,ac}$ with $\eta=10^{-3}$ from time step 2000.

\Cref{fig:L12_3d_planewave_cube}(b)--(d) shows the pulses crossing inside the cube at time step 3000 and, at time step 4500, the parts that passed through the cube lagging behind. The ranks stay at or below 16 ($\mathbf{E}$, $\mathbf{H}$) and 56 (auxiliary fields) in free space and reach 64 between time steps 2140 and 2560, as the pulses enter the cube (\Cref{fig:L12_3d_planewave_cube}(e)); the auxiliary fields, small except where the pulses cross the PML, reach it first, consistent with the noise-driven rank growth analyzed in \cite{zhou2025tensor}. Without the local smoothing, the simulation becomes unstable from about time step 4300 (\Cref{fig:L12_3d_planewave_cube}(f)). The simulation takes 10.1 hours (6.1~s per step at the maximum rank). QTT with the same settings agrees in field norm to within 0.5\% through time step 6000 and takes 1.7 times as long per step. After the pulses leave, the residual fields of the two simulations differ; without a reference at this size, we do not interpret them.}

%	
%	\begin{table}[tp!]
%		\footnotesize
%		\caption{Comparison of the HODLR and HODBF-enhanced eigen solver for the pillbox cavity model with $f$=20.472 GHz using 16 Cori Haswell nodes. }
%		\begin{center}
%			%\setlength\tabcolsep{1.5pt}
%			\begin{tabular}{|c|c|c|c|c|}
%				%\cline{2-6}
%				\hline			
%				Algorithm & HODLR & HODBF  & HODLR & Omega3P \\
%				\hline
%				Diameter & 13.6 $\lambda$ & 13.6 $\lambda$ & 13.6 $\lambda$ & 13.6 $\lambda$  \\
%				\hline				
%				\#of triangles & 213,792 & 213,792 & 53,240 & -  \\										
%				\hline
%				\#of tetrahedrons & - & - & - & 14,861,793  \\										
%				\hline	
%				Max edge length & 2.1 mm& 2.1 mm & 4 mm & 2.7 mm  \\			
%				\hline									
%				Basis & RWG & RWG & GWP($1$) & Nedelec\\				
%				\hline				
%				\#DoFs & 320,688 & 320,688 & 266,200 & 17,185,459 \\							
%				\hline				
%				Compression tolerance & $10^{-4}$ & $10^{-4}$ & $10^{-4}$ & - \\				
%				\hline				
%				Assembly\ylrev{\&compress} & 39 min &  8 min & 12 min &  16.8 min \\				
%				\hline						
%				Max. rank & 4113 &  380 & 4037 &  - \\		
%				\hline								
%				Inversion time & 54 s & 38 min & 52 s & -  \\			
%				\hline
%				Total memory  & 250 GB & 51.1 GB & 200.8 GB & 1.56 TB  \\
%				\hline				
%				Eigen solver time & 13.4 s & 49.2 s  & 11.9 s & 7.4 min   \\
%				\hline
%				\#of detected modes & 5 & 5  & 5 & 5   \\
%				\hline									
%			\end{tabular}
%		\end{center}
%		%\caption{Technical data for the examples in \cref{sec:acc}. $\lambda_0$ is the wavelength in the background medium.}
%		\label{tab:tech_data}
%		%\vspace{-20pt}
%	\end{table}
%	

	\section{Conclusion}\label{sec:con}
	This paper presents CTN-augmented 2D and 3D FDTD algorithms with PML boundaries\rev{, with local and global smoothing to keep the CTN ranks low. When the ranks stay bounded, as for plane-wave-like fields, the memory and the cost per time step grow logarithmically with $N$, plus $O(R^3)$ and $O(R^4)$ spine terms in 3D; the ranks grow slowly under refinement and roughly as $N$ for a point source at a fixed PPW. For a 3D dielectric sphere, CTN-FDTD agrees with standard FDTD to 2--4\% and needs about a third of the time per step of a QTT chain at comparable accuracy, and a $4096^3$ cube scatterer took 10 hours on one node, which is far more memory-efficient than standard FDTD. The algorithm also applies to other PMLs and, through polarization variables updated by the same operations, to dispersive media. Future work includes parallelization and more complex scatterers and media.}

	\section*{Acknowledgment}\label{sec:ack}
	This work was supported in part by the U.S. Department of Energy, Office of Science, Office of Advanced Scientific Computing Research, Applied Mathematics Competitive Portfolios program under Contract No. DE-AC02-05CH11231 at Lawrence Berkeley National Laboratory. This research was performed under an appointment to the Minority Serving Institutions Internship Program (MSIIP) administered by the Oak Ridge Institute for Science and Education (ORISE) for the National Nuclear Security Administration (NNSA) and U.S. Department of Energy (DOE). ORISE is managed by Oak Ridge Associated Universities (ORAU). All opinions expressed in this publication are those of the authors and do not necessarily reflect the policies and views of NNSA, DOE, ORISE or ORAU.

	\bibliographystyle{unsrt}
	\bibliography{References.bib}

\end{document}